\documentclass{proc-l}

\usepackage{amsmath,amssymb,amsthm,amscd}
\newtheorem{theorem}{Theorem}[section]

\theoremstyle{definition}

\theoremstyle{remark}
\newtheorem{remark}[theorem]{Remark}

\newcommand*\diff{\mathop{}\!\mathrm{d}}
\numberwithin{equation}{section}

\begin{document}

\title{Rigidity of entire self-shrinking solutions to K\"{a}hler-Ricci flow on complex plane}


%    author one information
\author{Wenlong Wang}
\address{School of Mathematical Sciences\\
Peking University\\
Science Building in Peking University, No 5. Yiheyuan Road, Beijing, P.R.China 100871}
\email{wwlpkumath@yahoo.com}
\thanks{The author is partially supported by CSC(China Scholarship Council)}

%    \subjclass is required.
\date{}

%    "Communicated by" -- provide editor's name; required.
\commby{Lei Ni}

%    Abstract is required.
\begin{abstract}
We show that every entire self-shrinking solution on $\mathbb{C}^1$ to the K\"{a}hler-Ricci flow must be generated from a
quadratic potential.
\end{abstract}

\maketitle

\section{Introduction}

In this short note, we prove the following result.
\begin{theorem}\label{thm}

Suppose that $u(x)$ is an entire smooth subharmonic solution  on $\mathbb{R}^{n}$ to the equation
\begin{equation}\label{LL}
\ln\Delta u=\frac{1}{2}x\cdot Du-u,
\end{equation}
then $u$ is quadratic.
\end{theorem}

For $n=2$, up to an additive constant, equation \eqref{LL} is equivalent to the one-dimensional case of the complex Monge-Amp\`ere equation 
\begin{equation}\label{KR}
\ln\det u_{i\bar{j}}=\frac{1}{2}x\cdot Du-u
\end{equation}
on $\mathbb{C}^m$. Any entire solution to \eqref{KR} leads to an entire self-shrinking solution 
\begin{equation}
v(x,t)=-tu\left(\frac{x}{\sqrt{-t}}\right)  \nonumber
\end{equation}
to a parabolic complex Monge-Amp\`ere equation 
\begin{equation}
v_t=\ln \det\left(v_{i\bar{j}}\right)\nonumber
\end{equation}
on $\mathbb{C}^{m} \times (-\infty,0)$, where $z^{i}=x^{i}+\sqrt{-1}x^{m+i}$. Note that above equation of $v$ is the potential equation of the K\"{a}hler-Ricci flow $\displaystyle\partial_tg_{i\bar{j}}=-R_{i\bar j}$. In fact, the corresponding metric $\left(u_{i\bar{j}}\right)$ is a shrinking K\"{a}her-Ricci (non-gradient) soliton.

Assuming a certain decay of $\Delta u$--a specific completeness condition, Q. Ding and Y.L. Xin have proved Theorem \eqref{thm} in \cite{DX}. Under the condition that the K\"{a}hler metric $\left(u_{i\bar{j}}\right)$ is complete, rigidity theorem for equation \eqref{KR} has been obtained by  G. Drugan, P. Lu and Y. Yuan  in \cite{DLY}. Similar rigidity results for self-shrinking solutions to Lagrangian mean curvature flows in pseudo-Euclidean space were obtained in \cite{CCY}, \cite{DX}, \cite{H} and \cite{HW}.

Our contribution is removing extra assumptions for the rigidity of equation \eqref{LL}. As in \cite{DLY} and \cite{DX}, the idea of our argument is still to prove the phase--$\ln\Delta u$ is constant. Then the homogeneity of the self-similar term on the right-hand side of equation \eqref{LL} leads to the quadratic conclusion. However, it's hard to construct a barrier function as in \cite{DLY} or to find a suitable integral factor as in \cite{DX} without completeness assumption. Taking advantage of the conformality of the linearized equation \eqref{LL}, we establish a second order ordinary differential inequality for $M(r)=\max_{|x|=r}\ln\Delta u(x)$ in the sense of comparison function. Then we prove that $M(r)$ blows up in finite time by Osgood's criterion unless $\ln\Delta u$ is constant.  

\section{proof}
\begin{proof}

Define the phase by
\begin{equation}
\phi(x)=\frac{1}{2}x\cdot Du(x)-u(x). \nonumber
\end{equation}
Taking two derivatives and using \eqref{LL}, we have 
\begin{equation}\label{phase}
\Delta\phi=\frac{e^\phi}{2}x\cdot D\phi.
\end{equation}

Define $M(r):[0,+\infty)\rightarrow\mathbb{R}$ by
\begin{equation}
M\left(r\right)=\max_{|x|=r}\phi\left(x\right). \nonumber
\end{equation}
Assuming $\phi\left(x\right)$ is not a constant, we prove that $M\left(r\right)$ blows up in finite time.

Since $M$ is locally Lipschitz, it is differentiable a.e. in $[0,+\infty)$. For all $r>0$, there exists a corresponding angle $\theta_r\in\mathbb{S}^{n-1}$ satisfying
\begin{equation}\label{M1}
M\left(r\right)=\phi\left(r,\theta_r \right).
\end{equation}

For $r'>0$ small enough, we have $M(r+r')\geq\phi(r+r',\theta _r)$ and $M(r-r')\geq \phi(r-r',\theta _r)$. It follows that
\begin{equation*}
\frac{M(r+r')-M(r)}{r'}\geq\frac{\phi(r+r',\theta_r)-\phi(r,\theta_r)}{r'}
\end{equation*}
and
\begin{equation*}
\frac{M(r)-M(r-r')}{r'}\leq\frac{\phi(r,\theta_r)-\phi(r-r',\theta_r)}{r'}.
\end{equation*}
Letting $r'\rightarrow 0$ in above two equations, we have
\begin{equation*}
\varlimsup M'_-(r)\leq \frac{\partial\phi}{\partial r}\left(r,\theta_r\right)\leq \varliminf M'_+(r).
\end{equation*}
So if $r>0$ is a differential point of $M$, we have
\begin{equation}\label{R2}
M'(r)=\frac{\partial\phi}{\partial r}\left(r,\theta_r\right).
\end{equation}
Because of the maximality of $\phi\left(r,\theta_r\right)$ among $\theta\in\mathbb{S}^{n-1}$, we have $\Delta_{\mathbb{S}^{n-1}}\phi(r,\theta_r)\leq 0$. Plugin this inequality into \eqref{phase}, we obtain
\begin{equation}\label{5}
\frac{\partial^2\phi}{\partial r^2}\left(r,\theta_r\right)+\frac{n-1}{r}\frac{\partial\phi}{\partial r}\left(r,\theta_r\right)\geq\frac{r}{2}\exp\left[\phi\left(r,\theta_r\right)\right]\cdot\frac{\partial\phi}{\partial r}\left(r,\theta_r\right).
\end{equation}

Fixing a positive $R_0$, for any $r\in[0, R_0]$, $\theta\in\mathbb{S}^{n-1}$, $t\in\left[0,1\right]$, we have the following Taylor's expansion
\begin{equation}
\phi\left(r+t,\theta\right)-\phi\left(r,\theta\right)\geq\frac{\partial\phi}{\partial r}\left(r,\theta_r\right)\cdot t+\frac{1}{2}\frac{\partial^2\phi}{\partial r^2}\left(r,\theta_r\right)\cdot t^2-Ct^3.\nonumber
\end{equation}
Here $C$ is a constant depends only on $R_0$, in fact we can choose $C=\max_{|x|\leq R_0+1}|D^3\phi(x)|$. 

We evaluate above inequality at $(r, \theta_r)$, where $r$ is a  differentiable point and $\theta_r$ is the corresponding critical angle. Using $M\left(r+t\right)\geq\phi\left(r+t,\theta_r\right)$, \eqref{M1}, \eqref{R2} and \eqref{5}, we have the following inequality that only involves $M$, namely
\begin{equation}
M\left(r+t\right)-M\left(r\right)\geq M'\left(r\right)t+\frac{1}{4}\left\{\exp\left[M\left(r\right)\right]r-\frac{2(n-1)}{r}\right\}M'\left(r\right)t^2-Ct^3.\nonumber
\end{equation}
Or equivalently,
\begin{equation}\label{M2}
\frac{1}{t^2}\left[M\left(r+t\right)-M\left(r\right)-M'\left(r\right)t\right]\geq\frac{1}{4}\left\{\exp\left[M\left(r\right)\right]r-\frac{2(n-1)}{r}\right\}M'\left(r\right)-Ct.
\end{equation}

Because inequality \eqref{M2} holds for $r\in\left[0, R_0\right]$ a.e., we can integrate it with respect to $r$ over any subinterval $[a,b]\subset[0,R_0]$, and get the following inequality for every $t\in[0,1]$,

\begin{align}\label{ieq}
&\frac{1}{t^2}\left\{\int_{b}^{b+t}M\left(r\right)\diff r-\int_{a}^{a+t}M\left(r\right)\diff r-\left[M\left(b\right)-M\left(a\right)\right]t\right\}\\\nonumber
\geq&\frac{1}{4}\int_a^b\left\{\exp\left[M\left(r\right)\right]r-\frac{2(n-1)}{r}\right\}M'\left(r\right)\diff r-C\left(b-a\right)t.
\end{align}

Choosing differentiable points $a$ and $b$ and letting $t \rightarrow0$ in \eqref{ieq}, we have
\begin{equation}\label{M3}
M'\left(b\right)-M'\left(a\right)\geq\frac{1}{2}\int_a^b\left\{\exp\left[M\left(r\right)\right]r-\frac{2(n-1)}{r}\right\}M'\left(r\right)\diff r.
\end{equation}

Since $R_0$ can be  arbitrarily large, in fact \eqref{M3} holds for all differentiable points $a, b\in\mathbb{R}_+$. 

We claim there exists $l_0>0$ such that $M'\left(r\right)>0$ at every differentiable point in $[l_0, +\infty)$. Otherwise, there exist an increasing sequence of differentiable points $\{r_k\}\subset\mathbb{R}_+$, and a sequence of corresponding critical angles $\{\theta_{k}\}\subset\mathbb{S}^{n-1}$ such that 
\begin{equation}
M'\left(r_k\right)=\frac{\partial\phi}{\partial r}\left(r_k,\theta_k\right)\leq 0, \,\,\,\, \mbox{and}\,\,\,\, \lim_{k\to\infty}r_k=+\infty.\nonumber
\end{equation}
Then according to Hopf's lemma, we know $\phi(x)$ is constant in $\displaystyle B_{r_k}\left(0\right)$. Since $r_k$ can be arbitrarily large, $\phi(x)$ is in fact a constant on the whole $\mathbb{R}^n$, which contradicts our assumption. 

So there exists a certain $l_0>0$, such that $\displaystyle M'\left(r\right)>0$ holds a.e. in $[l_0,+\infty)$. Then $M\left(r\right)$ monotonically increases on $[l_0,+\infty)$. When $a\geq l_1\triangleq l_0+n+2\exp\left[-M\left(l_0\right)\right]$, we have
\begin{align}\label{M4}
\int_a^b\left\{\exp\left[M\left(r\right)\right]r-\frac{2(n-1)}{r}\right\}M'\left(r\right)\diff r&>2\int_a^b \exp\left[M\left(r\right)\right]M'\left(r\right)\diff r\\\nonumber
&=2\left\{\exp\left[M\left(b\right)\right]-\exp\left[M\left(a\right)\right]\right\}.
\end{align}
Combining \eqref{M3} and \eqref{M4}, we obtain
\begin{equation}\label{M5}
M'\left(b\right)-M'\left(a\right)\geq\exp\left[M\left(b\right)\right]-\exp\left[M\left(a\right)\right].
\end{equation}

Above inequality holds for all differentiable points $a,b \in[l_1,+\infty)$.  
Choosing a differentiable point $l_2\geq l_1$, then $M'\left(r\right)\geq M'\left(l_2\right)>0$ holds a.e. in $\left[l_2,+\infty\right)$. Thus $M(r)\rightarrow+\infty$ as $r\rightarrow+\infty$. 

Then according to Osgood's criterion, $M(r)$ blows up in finite time, which contradicts the assumption that $\phi(x)$ is entire. Hence, we conclude $\phi(x)$ is constant.
Using $\phi(x)=\frac{1}{2}x\cdot Du(x)-u(x)$, we have 
\begin{equation}
\frac{1}{2}x\cdot D\left[u(x)+\phi(0)\right]=u(x)+\phi(0).\nonumber
\end{equation}

Finally, it follows from Euler's homogeneous function theorem that smooth $u(x)+\phi(0)$ is a homogeneous order 2 polynomial.
\end{proof}

\begin{remark}
From the proof, it's not hard to see that the theorem also holds for 
\begin{equation}
\Delta u=f\left(x\cdot Du-2u\right)\nonumber
\end{equation}
if $f\in C^1(\mathbb{R})$ is convex, monotone increasing, and $f^{-1}\in L^1([d,+\infty))$ for a certain $d\in\mathbb{R}$. Integrability condition for $f^{-1}$ is necessary. Otherwise, we have such counterexample: $f(x)\equiv x$ and
\begin{equation}
u(x)=\left(x_1^2-1\right)\int_0^{x_1}\frac{1}{s^2}(\exp\frac{s^2}{2}-1)\diff s-\frac{1}{x_1}(\exp\frac{x_1^2}{2}-1)-x_1.\nonumber
\end{equation}
\end{remark}
{\it Acknowledgement.} I would like to sincerely thank Professor Yu Yuan for suggesting this problem to me and for many heuristic discussions. The result was obtained when I was visiting the University of Washington. I would also like to thank  CSC(China Scholarship Council) for its support and the University of Washington for its hospitality. I am grateful to the referees for useful suggestions.

\bibliographystyle{amsplain}

\end{document}